\documentclass[ a4paper,
                reqno,
                twocolumn
                ]{article}
                
\usepackage{amsmath, amssymb, amsthm, mathtools}
\usepackage{array,booktabs}
\usepackage[skip=4pt,font=small,labelfont=bf,width=8cm]{caption}

\usepackage{extarrows}
\usepackage{tikz,hyperref}
\usetikzlibrary{matrix,arrows,decorations.markings,decorations.pathreplacing,petri,topaths}

\usepackage{sectsty}
\sectionfont{\large}

\newcommand{\letter}[1]{\mathtt{#1}}
\newcommand{\defgl}{\mathrel{\mathop:\!\!=}}
\newcommand{\IE}{\mathbb{E}}
\newcommand{\IN}{\mathbb{N}}
\newcommand{\IR}{\mathbb{R}}
\newcommand{\IX}{\mathbb{X}}
\newcommand{\IY}{\mathbb{Y}}
\newcommand{\IV}{\mathbb{V}}
\newcommand{\IZ}{\mathbb{Z}}
\newcommand{\CA}{\mathcal{A}}

\newcommand{\CD}{\mathcal{D}}
\newcommand{\CG}{\mathcal{G}}
\newcommand{\CH}{\mathcal{H}}
\newcommand{\CL}{\mathcal{L}}
\newcommand{\CN}{\mathcal{N}}
\newcommand{\CZ}{\mathcal{Z}}
\newcommand{\refpnt}{|}
\newcommand{\dif}{\mathop{}\!\mathrm{d}}
\newcommand{\ts}{\hspace{0.5pt}}
\newcommand{\II}{\ts\mathrm{i}\ts}

\DeclareMathOperator{\vol}{vol}
\DeclareMathOperator{\e}{e}

\newtheorem{thm}{Theorem}

\newtheorem{pro}[thm]{Proposition}

\begin{document}

\title{
\textbf{On a Family of Random Noble Means Substitutions}
}

\author{
Markus Moll
}

\date{
{\small Fakult\"at f\"ur Mathematik, Universit\"at Bielefeld \\ 
Universit\"atsstra{\ss}e 25, D--33615 Bielefeld, Germany}
}

\twocolumn[
\begin{@twocolumnfalse}
    \maketitle
    \begin{abstract} 
        \noindent In 1989, Godr{\`e}che and Luck \cite{luck} introduced the
        concept of local mixtures of primitive substitution rules along the
        example of  the well-known Fibonacci substitution and foreshadowed
        heuristic results on the topological entropy and the spectral type of
        the diffraction measure of associated point sets. In this contribution,
        we present a generalisation of this concept by regarding the so-called
        \textquoteleft noble means families', each consisting of finitely many
        primitive substitution rules that individually all define the same
        two-sided discrete dynamical hull.  We report about results in the
        randomised case on topological entropy, ergodicity of the two-sided
        discrete hull, and the spectral type of the diffraction measure of
        related point sets.
             
        \smallskip

        \bigskip
    \end{abstract}
\end{@twocolumnfalse}
]
{
  \renewcommand{\thefootnote}%
    {\fnsymbol{footnote}}
  \footnotetext[1]{email: \texttt{mmoll@math.uni-bielefeld.de}}
}

\bigskip

\section{The Setting}

Consider the binary alphabet $\CA_{2}^{} = \{ \letter{a},\letter{b} \}$. For
an arbitrary but fixed integer $m \ge 1$ and $0 \le i \le m$, we define the
\emph{noble means substitution (NMS)} rule $\zeta_{m,i}^{} \colon \CA_{2}^{}
\to \CA_{2}^{\ast}$ by
\begin{equation*}
    \zeta_{m,i}^{} \colon
    \left\{
    \begin{array}{lll}
        \letter{a}  & \mapsto & \letter{a}^{i} \letter{b} \letter{a}^{m - i},
        \\
        \letter{b}  & \mapsto & \letter{a},
    \end{array}
    \right.
    \quad \text{where} \quad
    M_{m}^{} \defgl
    \begin{pmatrix}
        m & 1 \\
        1 & 0
    \end{pmatrix}
\end{equation*}
is its (unimodular) substitution matrix. The family \[\CN_{m}^{} \defgl \bigl\{
\zeta_{m,i}^{} \mid m \in \IN, 0 \le i \le m \bigr\}\] is called a \emph{noble
means family} and each of its members is a primitive Pisot substitution with
inflation multiplier $\lambda_{m}^{} \defgl (m + \sqrt{m^{2} + 4})/2$ and
algebraic conjugate $\lambda_{m}^{\prime} \defgl (m - \sqrt{m^{2} + 4})/2$. The
two-sided discrete (symbolic) hull $\IX_{m,i}^{}$ of $\zeta_{m,i}^{}$ is
defined as the orbit closure of a fixed point in the local topology.  Each
$\IX_{m,i}^{}$ is reflection symmetric and aperiodic in the sense that it does
not contain any periodic element. For fixed $m \in \IN$, one observes that the
$\zeta_{m,i}^{}$ are pairwise conjugate and therefore all individual
$\IX_{m,i}^{}$ coincide; see \cite[Ch. 4]{bagrimm} for background.

Now, we fix $m \in \IN$ and a (strictly positive) probability vector
$\boldsymbol{p}_{m}^{} = (p_{0}^{},\ldots,p_{m}^{})$ and define a random
substitution $\zeta_{m}^{}$ on $\CA_{2}^{}$ by
\begin{equation*}
    \zeta_{m}^{} \colon 
    \left\{
    \begin{array}{lll}
        \letter{a} & \mapsto & 
        \left\{
        \begin{array}{cc}
            \zeta_{m,0}^{}(\letter{a}) & \text{with probability } p_{0}^{}, \\
            \vdots                     & \vdots                             \\
            \zeta_{m,m}^{}(\letter{a}) & \text{with probability } p_{m}^{},
        \end{array}
        \right. \\
        \letter{b} & \mapsto & \letter{a}.
    \end{array}
    \right.
\end{equation*}
We refer to $\zeta_{m}^{}$ as a \emph{random noble means substitution (RNMS)}.
Both the substitution matrix and the inflation multiplier are the same as in
the NMS case. We aim at the \emph{local mixture} of all members of
$\CN_{m}^{}$, which means that we independently apply $\zeta_{m}^{}$ to each
letter of some word $w \in \CA_{2}^{\IZ}$. In this case, the two-sided discrete
stochastic hull $\IX_{m}^{}$ is defined as the smallest closed and
shift-invariant subset of $\CA_{2}^{\IZ}$ with the property that $X_{m}^{}
\subset \IX_{m}^{}$, where
\begin{align*}
    X_{m}^{} \defgl 
    \Bigl\{
        w \in \mathcal{A}_{2}^{\mathbb{Z}} \mid & \;w \text{ is an accumulation }
        \\ &\text{ point of } 
    \bigl( 
        \zeta_{m}^{k}(\letter{a} \refpnt \letter{a}) 
    \bigr)_{k \in \IN_{0}} 
    \Bigr\}.
\end{align*}
Both $\IX_{m,i}^{}$ and $\IX_{m}^{}$ are completely characterised by the legal
subwords. Here, a word $w \in \CA_{2}^{\ast}$ is $\zeta_{m}^{}$-legal if there
is a $k \in \IN$ such that $w$ is a subword of at least one realisation of the
\emph{random variable} $\zeta_{m}^{k}(\letter{b})$. The set of
$\zeta_{m}^{}$-legal words of length $\ell$ is henceforth denoted by
$\CD_{m,\ell}^{}$. One can show that $\IX_{m,i}^{} \subsetneq \IX_{m}^{}$ by
considering the subword $\letter{bb}$ and that the system $(\IX_{m}^{},S)$,
where $S$ denotes the shift, is topologically transitive but not minimal.  Note
that $\IX_{m}^{}$ is invariant under alterations of $\boldsymbol{p}_{m}^{}$ as
long as $\boldsymbol{p}_{m}^{}$ is strictly positive.

\section{Topological entropy} \label{sec:entropy}

For $m \in \IN$ and $n \ge 3$, the set of \emph{exact RNMS words} is given by
\begin{equation} \label{equ:concat_rule}
    \CG_{m,n}^{} \defgl \bigcup_{i=0}^{m} \prod_{j=0}^{m}
    \CG_{m,n-1-\delta_{ij}^{}}^{},
\end{equation}
where $\CG_{m,1}^{} \defgl \{\letter{b}\}$, $\CG_{m,2}^{} \defgl
\{\letter{a}\}$ and $\delta_{ij}^{}$ denotes the Kronecker function. The
product in Eq. $\eqref{equ:concat_rule}$ is understood via concatenation of
words. Now, assume that $\boldsymbol{p}_{m}^{}$ is strictly positive. The
complexity function 
\begin{equation*}
    C_{m}^{} \colon \IN \to \IN, \quad
    \ell \mapsto \lvert \CD_{m,\ell}^{} \rvert
\end{equation*}
of $\zeta_{m}^{}$ is unknown, but the knowledge of the exact RNMS words is
enough \cite{nilsson} to compute the \emph{topological entropy} $\CH_{m}^{}$ of
$\zeta_{m}^{}$ for any $m \in \IN$ to be
\begin{align*}
    \CH_{m}^{} &= \lim_{n \to \infty} \frac{\log\bigl( C_{m}^{}(\ell_{m,n}^{})
    \bigr)}{\ell_{m,n}^{}} = \lim_{n \to \infty} \frac{\log\bigl( \lvert
    \CG_{m,n}^{} \rvert \bigr)}{\ell_{m,n}^{}} \\
    &= \frac{\lambda_{m}^{} - 1}{1 - \lambda_{m}^{\prime}} \sum_{i =
    2}^{\infty} \frac{\log\bigl( m(i-1) + 1 \bigr)}{\lambda_{m}^{i}} > 0,
\end{align*}
where $\ell_{m,n}^{}$ is the length of any word $w \in \CG_{m,n}^{}$. The
numerical values of $\CH_{m}^{}$ for $1 \le m \le 4$ are shown in Table
\ref{tab:entropy}.
\begin{center}
    \captionof{table}{Numerical values of $\CH_{m}^{}$ for $1 \le m \le 4$.}
    \begin{tabular}{ccccc}
        \specialrule{1pt}{0pt}{4pt}
        $m$ & $1$ & $2$ & $3$ & $4$ \\
        $\CH_{m}^{}$ & $0.44439$ & $0.40855$ & $0.37139$ &$0.33862$ \\
        \specialrule{1pt}{4pt}{1pt}
    \end{tabular}
    \label{tab:entropy}
\end{center}
One can prove that $\CH_{m}^{} > \CH_{m+1}^{}$ for all $m \in \IN$ and
$\mathcal{H}_{m}^{} \xrightarrow{m \to \infty} 0$.

\section{Ergodicity of $(\IX_{m}^{},S)$}

The known concept of the \emph{induced substitution} \cite[Ch. 5]{queff} that
acts on the alphabet of legal subwords of a fixed length can be generalised to
the stochastic setting of the RNMS case. One obtains a random substitution rule
$(\zeta_{m}^{})_{\ell}^{} \colon \CD_{m,\ell}^{} \to \CD_{m,\ell}^{\ast}$ and
one can prove that the induced substitution matrix $M_{m,\ell}^{}$ is a
primitive matrix which enables the application of Perron--Frobenius (PF)
theory.

\begin{figure}[t]
    \centering
    \begin{tikzpicture}[description/.style={fill=white,inner sep=7pt}]
    \matrix (D) [matrix of math nodes, row sep=2em,column sep=4.5em, text
    height=1.5ex, text depth=0.2ex,minimum width=2em,nodes in empty 
    cells,scale=0.45]
    { 
        \makebox[1cm]{$\IR$} & \makebox[1cm]{$\IR \times \IR$} &
        \makebox[1cm]{$\IR$}              \\
        \makebox[1cm]{$\IZ[\lambda_{m}^{}]$} & 
        \makebox[1cm]{$\mathcal{L}_{m}^{}$}  &
        \makebox[1cm]{$\IZ[\lambda_{m}^{}]$} \\ 
        \makebox[1cm]{$L$} &  & \makebox[1cm]{$L^{\star}$} \\ };
    \path[-stealth]
    (D-1-2) edge node[above] {$\pi_{1}^{}$} (D-1-1)
    (D-1-2) edge node[above] {$\pi_{2}^{}$} (D-1-3)
    (D-2-2) edge node[above] {$1-1$} (D-2-1)
    (D-2-2) edge node[above] {$1-1$} (D-2-3)
    (D-3-1) edge node[above] {$\star$} (D-3-3);
    \path[-] 
    (D-3-1) edge[double] (D-2-1)
    (D-3-3) edge[double] (D-2-3);
    \path[] 
    (D-2-1) edge node[description] {\rotatebox{90}{$\subset$}} (D-1-1)
    (D-2-2) edge node[description] {\rotatebox{90}{$\subset$}} (D-1-2)
    (D-2-3) edge node[description] {\rotatebox{90}{$\subset$}} (D-1-3);
    \node (1) at ( 3.5,0.65) {{\small dense}};
    \node (2) at (-3.5,0.65) {{\small dense}};
    \end{tikzpicture}
    \caption{Cut and project scheme for the noble means sets $\Lambda_{m,i}^{}$.}
    \label{fig:nms_cps}
\end{figure}
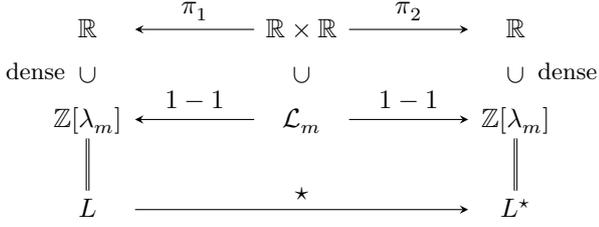

For fixed $m$, $\ell \in \IN$, let $w \in \CD_{m,\ell}^{}$ be a
$\zeta_{m}^{}$-legal word.  We define a \emph{shift-invariant probability measure}
$\mu_{m}^{}$ on the cylinder sets $\CZ_{k}^{}(w) = \{ v \in \IX_{m}^{} \mid
v_{[k,k + \ell - 1]}^{} = w \}$ for any $k \in \IZ$ by
\begin{equation} \label{equ:prob_meas}
    \mu_{m}^{} 
    \bigl( 
        \CZ_{k}^{}(w)
    \bigr) \defgl 
    \boldsymbol{R}_{m,\ell}^{}(w),
\end{equation}
where $\boldsymbol{R}_{m,\ell}^{}(w)$ is the entry of the (statistically
normalised) right PF eigenvector of $M_{m,\ell}^{}$ according to the word $w$.
\begin{thm}[{\cite{moll,bamoll}}] \label{thm:ergodic_discr}
    Let $\IX_{m}^{} \subset \CA_{2}^{\IZ}$ be the hull of the random noble
    means substitution for $m \in \IN$ and $\mu_{m}^{}$ the shift-invariant
    probability measure of Eq. $\eqref{equ:prob_meas}$ on $\IX_{m}^{}$.  For
    any $f \in L^{1}(\IX_{m}^{},\mu_{m}^{})$ and for an arbitrary but fixed $s
    \in \IZ$,
    \begin{equation*}
        \lim_{N \to \infty} \frac{1}{N} \sum_{i = s}^{N + s - 1} f(S^{i}x) =
        \int_{\IX_{m}^{}}f \dif \mu_{m}^{}
    \end{equation*}
    holds for $\mu_{m}^{}$-almost every $x \in \IX_{m}^{}$.
\end{thm}
Theorem \ref{thm:ergodic_discr} implies that $\mu_{m}^{}$ is ergodic. The proof
can be accomplished via an application of Etemadi's formulation of the strong
law of large numbers \cite{etemadi} and a suitable reorganisation of the
summation over the characteristic function of some cylinder set.

\section{Cut and project} \label{sec:cap}

The geometric realisations $\Lambda_{m,i}^{}$ of fixed points of each noble
means substitution can be derived as regular model sets within the cut and
project scheme $(\IR,\IR,\CL_{m}^{})$, where $\CL_{m}^{} \defgl \{
(x,x^{\prime}) \mid x \in \IZ[\lambda_{m}^{}] \}$; see Figure
\ref{fig:nms_cps}.  Here, the letters $\letter{a}$ and $\letter{b}$ are
identified with closed intervals of length $\lambda_{m}^{}$ and $1$,
respectively, and the left endpoints are chosen as control points.
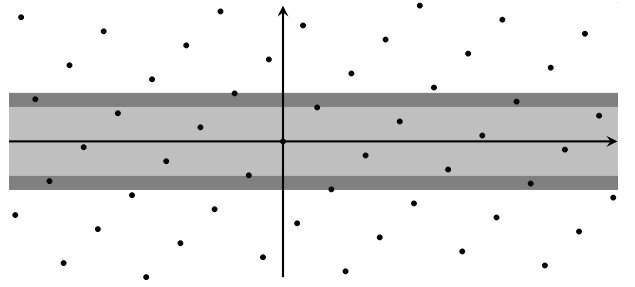
\begin{figure}[t]
    \centering
    \begin{tikzpicture}[scale=0.45]
    \tikzstyle{prj}=[color=white,fill=#1,line width=1.1pt] 
    \def\InfMult{2.4142}   
    \def\InfMultC{-0.4142} 
    \def\lSuperBound{-1.4142}
    \def\uSuperBound{1.4142}
    \def\Rad{0.1cm} 
    \def\FillUpper{red!50} 
    \def\FillLower{gray!50} 
    \def\FillSuper{black!50}
    \clip (-8,-4.1) rectangle (9.78,4.1); 
    \draw[color=\FillSuper,fill=\FillSuper] (-10,\lSuperBound) rectangle
    (10,\uSuperBound);
    \draw[color=\FillLower,fill=\FillLower] (-10,-1) rectangle
    (10,1);
    \draw[-stealth,thick] (-10, 0) -- (9.78,0);
    \draw[-stealth,thick] (0,  -4) -- (0,   4);
    \foreach \x in {-5,...,5} {
        \foreach \y in {-5,...,5} {
            \filldraw (\x+\InfMult*\y,\x+\InfMultC*\y) circle (0.07cm);
        }
    }
    \end{tikzpicture}
    \caption{The strip $\IR \times
        \bigcup_{i=0}^{m} W_{m,i}^{}$ (light) is strictly included in
        $\IR \times W_{m}^{}$ (dark). Here, this is illustrated for
        $m = 2$.}
    \label{fig:superwindow}
\end{figure}
The windows $W_{m,i}^{}$ for $\Lambda_{m,i}^{}$, in the generic cases $0 < i <
m$, are
\begin{equation*}
    W_{m,i}^{} \defgl i\tau_{m}^{} + [\lambda_{m}^{\prime},1] \quad\text{with}\quad
    \tau_{m}^{} \defgl -\frac{1}{m}(\lambda_{m}^{\prime} + 1),
\end{equation*}
while in the singular cases $i=0$ and $i=m$, one finds
\begin{alignat*}{3}
    W_{m,0}^{(\letter{a}|\letter{a})} &\defgl [\lambda_{m}^{\prime},1) , &
        W_{m,0}^{(\letter{a}|\letter{b})}
    &\defgl (\lambda_{m}^{\prime},1], \label{equ:model:1} \\
    \quad W_{m,m\vphantom{0}}^{(\letter{a}|\letter{a})} &\defgl
    (-1,-\lambda_{m}^{\prime}], & \quad
    W_{m,m\vphantom{0}}^{(\letter{b}|\letter{a})} &\defgl
    [-1,-\lambda_{m}^{\prime}),
\end{alignat*}
distinguished according to the legal two-letter seeds. Now, one can prove that
each member of the continuous RNMS hull $\IY_{m}^{}$ is a relatively dense
subset of an element of the LI class of the model set $\Theta(W_{m}^{})$, within
the cut and project scheme of Figure \ref{fig:nms_cps}, with window $W_{m}^{} =
[\lambda_{m}^{\prime} - 1, 1 - \lambda_{m}^{\prime}]$ and therefore a Meyer set
by \cite[Thm. 9.1]{moody}. The volume of the interval $W_{m}^{}$ is minimal with
this property, and it strictly contains $\bigcup_{i=0}^{m} W_{m,i}^{}$; see
Figure \ref{fig:superwindow} for an illustration in the case of $m=2$.

Consequently, each geometric realisation of a random noble means word is a 
naturally arising instance of a Meyer set with entropy.

\section{Diffraction}

The diffraction of the NMS cases is well understood due to their
characterisation as regular model sets \cite[Ch. 9]{bagrimm} whereas the results
presented in \cite{strungaru,baleri}, suggest the presence of a continuous part
in the diffraction spectrum in the RNMS case.

\begin{figure}[t]
    \centering 
    \includegraphics[width=8cm]{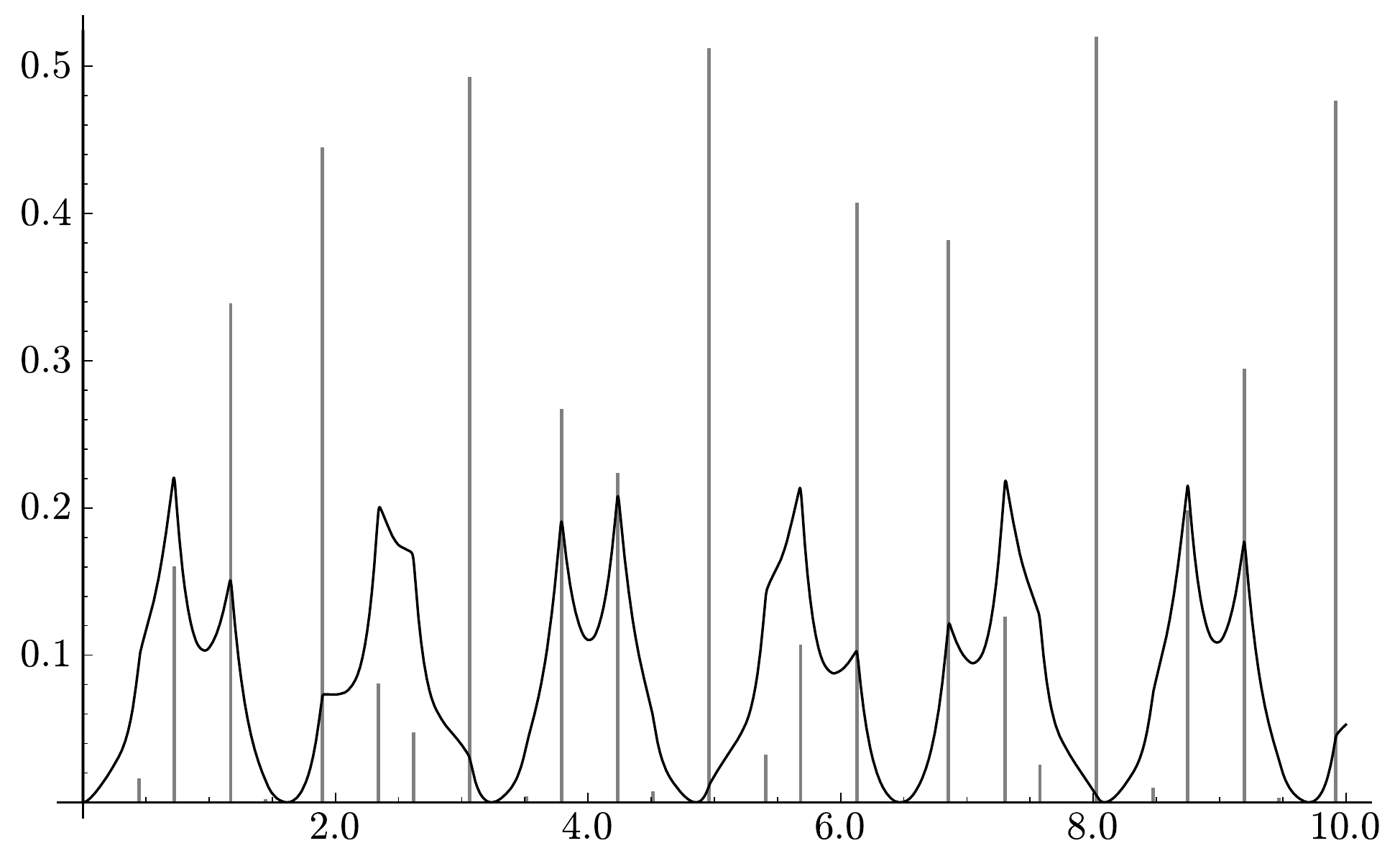}
    \caption{The pure point (light) and absolutely continuous (dark) part of 
        $\widehat{\gamma_{1}^{}}$ with $\boldsymbol{p}_{1}^{} = 
        (1/2,1/2)$ is illustrated.}
    \label{fig:abscont_pp}
\end{figure}

Because of Theorem \ref{thm:ergodic_discr}, the suspension \cite{cornfeld,ward}
$\nu_{m}^{}$ of $\mu_{m}^{}$ on $\IY_{m}^{}$ leads to a continuous and ergodic
dynamical system $(\IY_{m}^{},\IR,\nu_{m}^{})$. Now, let $\delta_{\Lambda}^{}
\defgl \sum_{x \in \Lambda} \delta_{x}^{}$ be the Dirac comb for a random noble
means set $\Lambda \in \IY_{m}^{}$. One can compute the autocorrelation of
$\delta_{\Lambda}^{}$ to be $\nu_{m}^{}$-almost surely given by
\begin{equation*}
    \gamma \defgl \lim_{R \to \infty} \frac{\delta_{\Lambda_{R}}^{} \ast 
    \widetilde{\delta_{\Lambda_{R}}^{}}}{\vol(B_{R}^{})} =
    \IE(\delta_{\Lambda}^{} \circledast 
    \widetilde{\delta_{\Lambda}^{}}),
\end{equation*}
where $\Lambda_{R}^{} \defgl \Lambda \cap B_{R}^{}(0)$ and $\circledast$
denotes the volume-averaged convolution by balls.  The diffraction measure is
given by the Fourier transform of $\gamma$ and reads
\begin{equation*}
    \widehat{\gamma} = \lim_{R \to \infty} \frac{1}{\vol(B_{R}^{})} \lvert
    \IE(\widehat{\delta_{\Lambda_{R}}}) \rvert^{2} + \lim_{R \to \infty}
    \frac{1}{\vol(B_{R}^{})} \IV(\widehat{\delta_{\Lambda_{R}}}),
\end{equation*}
where $\IE$ and $\IV$ refer to mean and variance with respect to the measure
$\nu_{m}^{}$. Now, the following two key properties finally lead to an explicit
expression for $\widehat{\gamma}$. 
\begin{enumerate}
    \item[\textbullet] It is enough to study $\widehat{\gamma}$ on the basis of
        exact RNMS words, as defined in Eq. $\eqref{equ:concat_rule}$, because
        $\zeta_{m}^{}$-legality of a word $w \in \CA_{2}^{\ast}$ means that $w$
        is a subword of a word in $\CG_{m,n}^{}$ for a suitably chosen $n \in
        \IN$.
    \item[\textbullet] It is not difficult to prove that 
        \begin{equation} \label{equ:stoch}
            \CG_{m,n}^{} = 
            \bigl\{ 
                w \in \CA_{2}^{\ast} \mid w = \zeta_{m}^{n-1}(\letter{b})
            \bigr\} 
        \end{equation} 
        and even more that the two stochastic processes, based on the
        substitution rule and the concatenation rule, are equal. Note that the
        equality in Eq. $\eqref{equ:stoch}$ means that there is at least one 
        realisation
        of the random variable $\zeta_{m}^{n-1}(\letter{b})$ that equals $w$.
\end{enumerate}
\begin{figure}[t]
    \centering 
    \includegraphics[width=8cm]{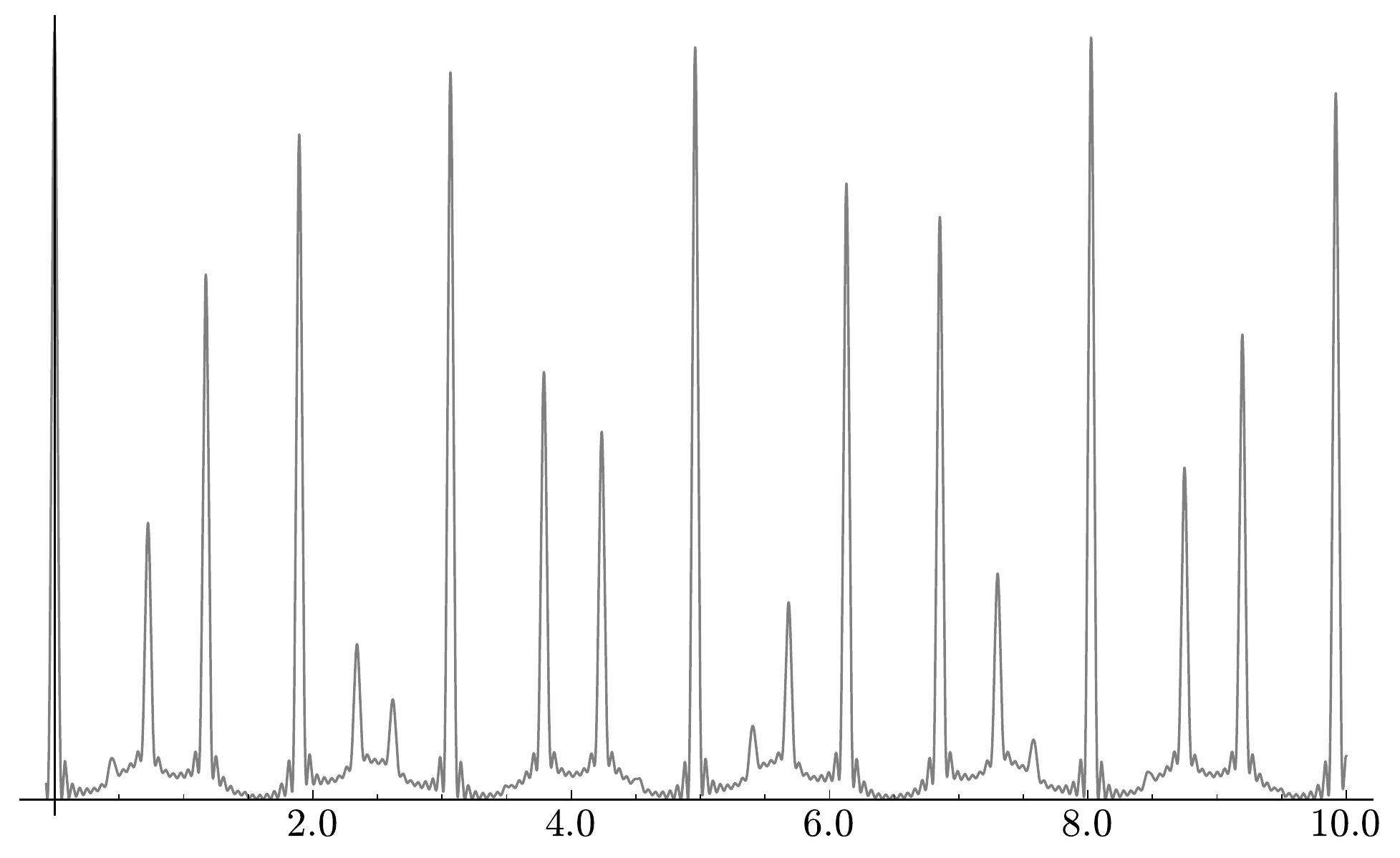}
    \caption{Approximation of the diffraction measure for $m=1$ and with $\boldsymbol{p}_{1}^{} = 
        (1/2,1/2)$ is illustrated.}
    \label{fig:full_diff}
\end{figure}
For convenience, we restrict to $m=1$ in the following and define for $n \ge 2$
the complex-valued random variable $X_{n}^{}(k)$ by
\begin{equation*}
    X_{n}^{}(k) \defgl
    \left\{
    \begin{array}{ll}
        X_{n-2}^{}(k) + \e^{-2\pi\II k \lambda_{1}^{n-2}} X_{n-1}^{}(k),
        & \langle p_{0}^{} \rangle, \\
        X_{n-1}^{}(k) + \e^{-2\pi\II k\lambda_{1}^{n-1}} X_{n-2}^{}(k),
        & \langle p_{1}^{} \rangle,
    \end{array}
    \right.
\end{equation*}
with $X_{0}^{}(k) = \e^{-2\pi\II k}$ and $X_{1}^{}(k) = \e^{-2\pi\II
k\lambda_{1}^{}}$. Here, $X_{n}^{}(k)$ corresponds to exact RNMS words in
$\CG_{1,n+1}^{}$. Therefore, we consider averaging over the sequence
$L_{n}^{} = \lambda_{1}^{n}$ and find the following result.
\begin{pro}[{\cite{moll}}]\label{pro:radon_niko_density}
    For any $n \in \IN$, consider the function $\phi_{n}^{} \colon \IR \to 
    \IR_{+}^{}$, defined by
    \begin{equation*}
        \phi_{n}^{}(k) \defgl \frac{1}{L_{n}^{}} \IV\bigl( X_{n}^{}(k) \bigr).
    \end{equation*}
    The sequence $(\phi_{n}^{})_{n \in \IN}^{}$ converges uniformly to the
    continuous function $\phi \colon \IR \to \IR_{+}^{}$, with
    \begin{equation}\label{equ:radon_niko_density}
        \phi(k) \defgl \frac{2p_{0}^{}p_{1}^{} \lambda_{1}}{\sqrt{5}} 
        \sum_{i=2}^{\infty} \lambda_{1}^{-i} \Psi_{i}^{}(k).
    \end{equation}
\end{pro}
Here, $\Psi_{n}^{} \colon \IR \to \IR_{+}^{}$ is a bounded and smooth function that 
monotonically decreases in $n$, defined by
\begin{equation*}
    \Psi_{n}^{}(k) \defgl
    \frac{1}{2}\bigl\lvert (1-\e_{n-2}^{}) \IE_{n-1}^{} -
    (1-\e_{n-1}^{}) \IE_{n-2}^{} \bigr\rvert^{2},
\end{equation*}
where $\e_{n}^{} \defgl \e^{-2\pi\II k\lambda_{1}^{n}}$ and $\IE_{n}^{} \defgl 
\IE(X_{n}^{}(k))$. This fixes the absolutely continuous part of $\widehat{\gamma}$. 
The pure point part can be computed via the recursion relation
\begin{equation} \label{equ:pp_rek}
    \IE_{n}^{}
    = (p_{1}^{} + p_{0}^{} \e_{n-2}^{})\IE_{n-1}^{} + (p_{0}^{} + p_{1}^{}
    \e_{n-1}^{})\IE_{n-2}^{},
\end{equation}
where $\IE_{0}^{} \defgl \e^{-2\pi\II k}$ and $\IE_{1}^{} \defgl \e^{-2\pi\II k 
\lambda_{1}^{}}$. This yields
\begin{equation*}
    \widehat{\gamma}(\{k\}) = \lim_{n \to \infty} \frac{1}{L_{n}^{2}} \lvert 
    \IE(X_{n}^{}(k)) \rvert^{2}
\end{equation*}
and an approximation of $(\widehat{\gamma})_{\mathrm{pp}}^{}$ and 
$(\widehat{\gamma})_{\mathrm{ac}}^{}$ is illustrated in Figure \ref{fig:abscont_pp} 
together with a sketch of the full diffraction, based on the recursion of Eq.
$\eqref{equ:pp_rek}$ with $n = 6$, in Figure \ref{fig:full_diff}.
Considering the Lebesgue decomposition $\widehat{\gamma} =
(\widehat{\gamma})_{\mathrm{pp}}^{} + (\widehat{\gamma})_{\mathrm{ac}}^{} +
(\widehat{\gamma})_{\mathrm{sc}}^{}$, we find that
\begin{equation*}
    \widehat{\gamma} = (\widehat{\gamma})_{\mathrm{pp}}^{} + \phi(k) \lambda,
\end{equation*}
where $\lambda$ denotes the Lebesgue measure and $\phi$ the density function of
Eq.  $\eqref{equ:radon_niko_density}$.

It is possible to compute the pure point part from the recursion relation in 
Eq. $\eqref{equ:pp_rek}$. Another interesting approach comes
from the theory of iterated function systems and inflation-invariant measures.
Here, one finds that
\begin{equation*}
    (\widehat{\gamma})_{\mathrm{pp}}^{} = \sum_{k \in
        \mathcal{L}_{1}^{\circledast}} \lvert
        \widehat{\eta_{\letter{a}}^{}}(-k^{\prime}) +
        \widehat{\eta_{\letter{b}}^{}}(-k^{\prime}) \rvert^{2} \delta_{k}^{},
\end{equation*}
where $\mathcal{L}_{1}^{\circledast} = \pi_{1}^{}(\mathcal{L}_{1}^{\ast}) =
\IZ[\lambda_{1}^{}]/\sqrt{5}$, with $\mathcal{L}_{1}^{\ast}$ the dual
lattice of $\mathcal{L}_{1}^{}$, is the Fourier module.  In the following, we
write $\xi \defgl \lambda_{1}^{\prime}$. The invariant measures
$\widehat{\eta_{\letter{a}}^{}}$, $\widehat{\eta_{\letter{b}}^{}}$ can be
approximated via the recursion relation
\begin{equation*}
    \begin{pmatrix}
        \widehat{\eta_{\letter{a}}^{}}(k) \\
        \widehat{\eta_{\letter{b}}^{}}(k)
    \end{pmatrix}
    =
    \lvert \xi \rvert^{n}
    \Bigl(
        \prod\limits_{\ell=1}^{n} p_{0}^{} A_{\ell}^{}(k) + p_{1}^{}
        B_{\ell}^{}(k)
    \Bigr)
    \begin{pmatrix}
        \widehat{\eta_{\letter{a}}^{}}(k\xi^{n}) \\
        \widehat{\eta_{\letter{b}}^{}}(k\xi^{n})
    \end{pmatrix},
\end{equation*}
where the matrices $A_{\ell}^{}(k)$ and $B_{\ell}^{}(k)$ are given by
\begin{equation*}
    \begin{pmatrix*}[l]
        \mathrm{e}^{-2\pi\mathrm{i}k\xi^{\ell-1}} & 1 \\
        1                                         & 0
    \end{pmatrix*}
    \quad \text{and} \quad
    \begin{pmatrix*}[l]
        1                                       & 1 \\
        \mathrm{e}^{-2\pi\mathrm{i}k\xi^{\ell}} & 0
    \end{pmatrix*}.
\end{equation*}
As $\xi^{n} \to 0$ for $n \to \infty$, an appropriate choice of the eigenvector 
$\bigl( \widehat{\eta_{\letter{a}}^{}}(0), \widehat{\eta_{\letter{b}}^{}}(0) 
\bigr)^{T}$ for the 
equation
\begin{equation*}
    \begin{pmatrix*}
        1 & 1 \\
        1 & 0
    \end{pmatrix*}
    \begin{pmatrix*}
        \widehat{\eta_{\letter{a}}^{}}(0) \\
        \widehat{\eta_{\letter{b}}^{}}(0)
    \end{pmatrix*}
     =
     \lambda_{1}^{}
     \begin{pmatrix*}
         \widehat{\eta_{\letter{a}}^{}}(0) \\
         \widehat{\eta_{\letter{b}}^{}}(0)
     \end{pmatrix*}
\end{equation*}
fixes the recursion. Since $\widehat{\eta_{\letter{a}}^{}}(0) + 
\widehat{\eta_{\letter{b}}^{}}(0)$ must be the point density of some random golden 
means set, which always is $\lambda_{1}^{}/\sqrt{5}$, one finds 
$\widehat{\eta_{\letter{a}}^{}}(0) = 1/\sqrt{5}$ and 
$\widehat{\eta_{\letter{b}}^{}}(0) = (\lambda_{1}^{} - 1)/\sqrt{5}$.

\begin{figure}
    \centering 
    \includegraphics[width=8cm]{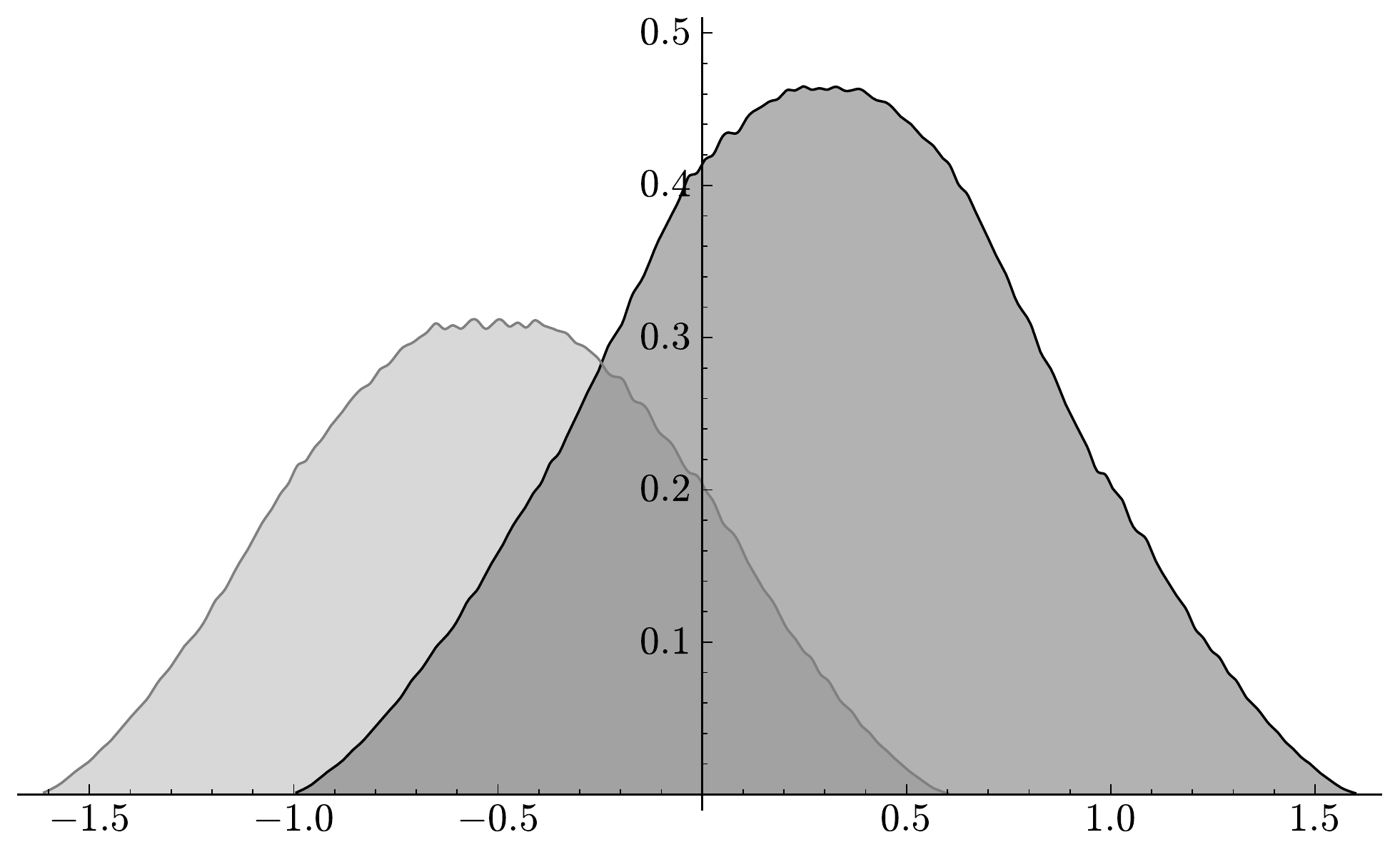}
    \caption{Distribution of control points generated by $\letter{a}$ (dark)
    and $\letter{b}$ (light) in the internal space in the case of
    $\boldsymbol{p}_{1}^{} = (1/2, 1/2)$. The plot is generated by a lift of 
    $\zeta_{1}^{32}(\letter{b})$ (i.e. 2178309 points) to the internal space.}
    \label{fig:distr_control}
\end{figure}

The distribution of control points in the internal space distinguished for 
$\letter{a}$ and $\letter{b}$, respectively, is illustrated in Figure 
\ref{fig:distr_control}.

\section*{Acknowledgement}

The author wishes to thank Michael Baake, Tobias Jakobi and Johan Nilsson for
helpful discussions. This work is supported by the German Research Foundation
(DFG) via the Collaborative Research Centre (CRC 701) through the faculty of
Mathematics, University of Bielefeld.

\end{document}